\documentclass[graybox]{svmult}

\usepackage{comment}
\usepackage{tikz}
\usepackage{type1cm}        
\usepackage{makeidx}         
\usepackage{graphicx}        
                            
\usepackage{multicol}        
\usepackage[bottom]{footmisc}

\usepackage{newtxtext}       
\usepackage{newtxmath}      

\usepackage{bm}
\if{false}
\fi
\newtheorem{rmk}[definition]{\bf Remark~\thermk}
\newtheorem{prp}[theorem]{\bf Proposition~\thetheorem}

\newcommand{\bea}{\begin{eqnarray}}
\newcommand{\eea}{\end{eqnarray}}
\newcommand{\bean}{\begin{eqnarray*}}
\newcommand{\eean}{\end{eqnarray*}}
\newcommand{\beq}{\begin{equation}}
\newcommand{\eeq}{\end{equation}}
\newcommand{\bac}{\begin{array}{c}}
\newcommand{\ball}{\begin{array}{ll}}
\newcommand{\ea}{\end{array}}

\renewcommand{\thetheorem}{\thesection.\arabic{theorem}}
\if{false}

\fi

\renewcommand{\thermk}{{ \thesection.\arabic{definition}}}

\def\a{{\bf a}}

\def\x{{\bf x}}

\def\y{{\bf y}}

\def\e{{\bf e}}
\def\1{{\bf 1}}
\def\0{{\bf 0}}

\def\RR{{\mathbb R}}
\def\RRk{{\RR^k}}

\if{false}

\renewcommand{\thermk}{{\rm \thesection.\arabic{rmk}}}

\fi

\def\al{\alpha}

\numberwithin{equation}{section}
\numberwithin{theorem}{section}

\def\({\left(} 
 
\def\){\right)}

\makeindex             
                       
\begin{document}

\title*{On Min-Max affine approximants of convex or concave real valued functions from $
\mathbb R^k$, Chebyshev equioscillation and graphics.}
\titlerunning{On Min-Max affine approximants}

\author{Steven B. Damelin, David L. Ragozin and Michael Werman}
\authorrunning{S. B. Damelin, D. L. Ragozin and M. Werman
}

\institute {Steven B. Damelin \at Department of Mathematics, University of Michigan, 530 Church Street, Ann Arbor, MI 48109, USA. 
 \email{damelin@umich.edu}\and David L. Ragozin \at Department of Mathematics, University of Washington, Seattle, WA 98195, USA. \email{rag@uw.edu}\and Michael Werman \at Department of Computer Science, The Hebrew University, 91904, Jerusalem, Israel. \email{michael.werman@mail.huji.ac.il}}

\maketitle
 
\abstract {We study  Min-Max affine approximants of a continuous convex or concave function
$f:\Delta\subset \mathbb R^k\xrightarrow{} \mathbb R$ where $\Delta$ is a convex compact subset of $\mathbb R^k$. 
In the case when $\Delta$ is a \emph{simplex} we prove that there is a \emph{vertical} translate of the supporting hyperplane in $\RR^{k+1}$ of the graph of $f$ at the vertices which is the unique best affine approximant to $f$ on $\Delta$.
For $k=1$, this result provides an extension of the Chebyshev equioscillation theorem 
for linear approximants. Our result has interesting connections to the computer graphics problem of rapid rendering of projective transformations.}
\medskip


\section{Introduction}
\setcounter{equation}{0}

We will work in $\RR^k , k\ge1$ where $\x \in \RRk$ is the column vector $[x_1\, x_2\, \cdots\, x_k]^{T}$ and $T$ denotes transpose. A function $g:\RRk \rightarrow \RR$ is an \emph{affine} function provided there exists ${\bm \alpha} \in \RRk,\,{\rm and}\, \beta \in \RR$ such that $g(\x)= {\bm \alpha}^T \x + \beta$.

\subsection{ Min-Max approximation}

In this paper, we are interested in  Min-Max approximation to a continuous $f:\Delta\subset \mathbb R^k\xrightarrow{} \mathbb R$ 
by affine approximants to $f$. See for example \cite{D,Da,L,Nik,Tem}. By imposing the restriction $\Delta$ is a simplex with non-empty interior, we obtain a complete characterization of these approximants and an explicit formula for a unique best approximant.
Even in the case of an interval, $k=1$, our main result Theorem \ref{t:simplexD} provides an extension of [\cite{Da}, Corollary 7.6.3].
For $k=1$, our main result also provides an extension of the Chebyshev equioscillation theorem for linear approximants with an explicit unique formula for the best approximant. See Section 4. 
We show interesting connections of Theorem \ref{t:simplexD} to graphics. See Section 5. 

In order to state our results, we need the following notation.
For a continuous function $g:\mathbb R^k\to \mathbb R$, and simplex $\Delta$, we adopt the usual convention of: 
\[
\left\|g(\x)
\right\|_{\infty (\Delta)}:= \displaystyle\sup_{\x\in \Delta}\left|g(\x)
\right|.
\]

We will answer the following
\subsection{Problem}

Let $f:\Delta\subset \mathbb R^k \to \mathbb R$ be a continuous function where $\Delta$ is a compact domain in $\mathbb R^k$. 
Find conditions on $\Delta$ and $f$ which allow for the construction of the Min-Max affine
approximation problem (\ref{e:minmaxmetric}) to $f$ over 
$\Delta$. 
Equivalently, find conditions on $\Delta$ and $f$ which allow for the explicit construction of an ${\bm\alpha} \in \RRk$ and a $\beta \in \RR$ which solve
\beq
\displaystyle
\min_{\left\{{\bm \al},\beta\right\}} \left\| f(\x) -({\bm \al}^T \x + \beta)
\right\|_{\infty (\Delta)}.
\label{e:minmaxmetric} 
\eeq

\section{Main result: Theorem \ref{t:simplexD}}
\setcounter{equation}{0}

Our {\bf main result} is

\begin{theorem}
Let $\{\a_1,\dots, \a_{k+1}\}$ be $k+1$ affinely independent points in $\mathbb R^k$ so that their convex hull $\Delta = {\rm CH}(\a_1 \dots \a_{k+1})$ is a $k$-simplex and assume 
$f:\Delta\subset \mathbb R^k \rightarrow \mathbb R$ is a continuous convex or concave function over $\Delta$. Then the  Min-Max affine approximant to $f$ over $\Delta$ is the hyperplane average $\sigma:=\frac{\pi+\rho}{2}$,
where $\pi$ is the affine hyperplane ${\rm AS}((\a_1, f(\a_1)) \dots (\a_{k+1},f(\a_{k+1})))$ in $\mathbb R^{k+1}$ and
$\rho$ is the supporting hyperplane to the graph of $f$ parallel to $\pi$.

Here 
${\rm AS}$ denotes 
the affine span and
any hyperplane $\tau$ is identified with its graph $\{(\x,\tau(\x)) \in \RR^{k+1}: x\in \RR^k \}$. 
\label{t:simplexD}
\end{theorem}
 
Theorem~\ref{t:simplexD} belongs to an interesting class of related optimization problems which can be found for example in
\cite{ANW, Br, Da, D, Korn, L, Mir, Nik, Tem}.
\medskip

We know of no convex domain other than a $k$-simplex where we can 
generate a hyperplane $\pi$ as in Theorem~\ref{t:simplexD} for which we can verify that the graph of $f$ over the domain lies entirely above or entirely below that hyperplane.
\medskip

We now present two remarks regarding
Theorem~\ref{t:simplexD}. 

\begin{rmk}
The following stronger version of Theorem~\ref{t:simplexD}
holds which has an identical proof as Theorem~\ref{t:simplexD}.
\end{rmk}

\begin{theorem}
Given a continuous function $f$ on a simplex $\Delta$ with the following property: The graph of $f$ lies entirely above or below its secant hyperplane $\pi$ through the graph points of $f$ over the vertices of $\Delta$. Then an expression for the  Min-Max affine approximation of $f$ on the simplex $\Delta$ has graph given by the hyperplane $\pi+d$ where 2d is the non-zero extremum value of $f-\pi$ on $\Delta$. 
% and where $\pi$ is the secant hyperplane to $f$ through the vertices of the simplex $\Delta$
An alternative definition of $\pi$ is the affine interpolant to $f$ at the $k+1$ vertices of the simplex $\Delta$.
\label{t:simplexDD}
\end{theorem}

\begin{rmk}

In this remark, we speak to a geometric description of the hyperplane average in Theorem~\ref{t:simplexD} and write formulae for the optimal ${\bm \al}$ and $\beta$ in the notation of Problem 1.2, expression  (1).

Let the vector ${\bm \al}'$ and the scalar $\beta'$ be defined as the solution to the following equation.

\begin{equation}
\begin{bmatrix}
\a_{1}^T &  1  \\
. &  1  \\
. &  1  \\
. &  1  \\
\a_{k+1}^T &  1
\end{bmatrix}
\begin{bmatrix}
{\bm \al}' \\

\beta'
\end{bmatrix}
=\begin{bmatrix}
f({\a_{1}}) \\
. \\
.  \\
. \\
f({\a_{k+1}})
\end{bmatrix}.
\end{equation}

Secondly define:

\begin{equation}
 \beta''=
    \begin{cases}
      \min_{{\bf x}\in \Delta}(f(\x) -({\bm \al}'^T \x+\beta')), & \text{if $f$ is convex}.\\
      \max_{{\bf x}\in \Delta}(f(\x) -({\bm \al}'^T \x+\beta')), & \text{if $f$ is concave}.
    \end{cases}       
\end{equation}

Then set ${\bm \al}: = {\bm \al}'$ and $\beta:=\frac{1}{2}(\beta'+\beta'')$. It is easily checked that Problem 1.2, expression (1) is optimized at 
$({\bm \al},\beta)$.

We end this remark by saying that the computation of the optimal $\beta$ above costs a minimization of a convex function (a linear shift of $\pm f$)
over the set $\Delta$.
\end{rmk}

\section{The proof of Theorem~\ref{t:simplexD}}
\setcounter{equation}{0}

The key ideas in our proof of Theorem~\ref{t:simplexD} 
are the two equivalent definitions of a convex function. $f:\Delta \subset \mathbb R^k\to \mathbb R$ is convex provided $\textnormal {for all } \gamma_i\geq 0: \Sigma_{i=1}^{k+1}\gamma_i =1,$ and any affine independent set of $k+1$ points $\{{\a}'_1,\ldots,{\a}'_{k+1}\}$, $f(\Sigma_{i=1}^{k+1}\gamma_i {\a}'_i)\leq \Sigma_{i=1}^{k+1}\gamma_if({\a}'_i)$, i.e. {the graph of $f$ over the convex hull of $k+1$ affinely independent points lies below the convex hull of the $k+1$ image points $\{[{\a}'_i{}^T\,\,f({\a}'_i)]^T: 1\leq i \leq k+1\}$, or equivalently, $\textnormal {for each } \x\in \Delta,$ $ \textnormal{and each support plane at }[\x^T\,\,f(\x)]^T$ $ \textnormal{ with slope } {\bf \eta} \in \RR^k $, $ f(\y)\geq f(\x)+{\bf \eta}^T(\y-\x),$ $ \textnormal{for all } \y\in \Delta$, i.e. the graph of $f$ lies on or above any supporting hyperplane.

For our proof of Theorem~\ref{t:simplexD}, we advise the reader to use Figure 1 for intuition.
\begin{proof}
First assume that $f$ is convex, so its graph over $\Delta$ lies on or below the hyperplane $\pi$. Let $-2d := \inf \{l: \pi + l\e_{k+1} \cap graph(f) \neq \varnothing\}$. 
Then any admissible $l$ must satisfy $l\leq 0$ and so this in turn means $d\geq 0$.
Here, $\e_{k+1}$ is the $k+1^{st}$ basis vector in $\mathbb R^{k+1}$ so each admissible $l$ gives a downward translate of $\pi$ with non-empty intersection with the graph of $f$. Since $\Delta$ is compact and $f$ is continuous, the $\inf$ is actually a $\min$ and so their exists at least one graph point $(y, f(y))$ on $\rho:=\pi - 2d$.

Thus $\rho$ is a supporting hyperplane (tangent plane if $f$ is differentiable at $\y$). 
Moreover, at the points $\a_i$ an easy computation shows $\pi - d = (2\pi - 2d)/2 = (\pi+\rho)/2= \rho+d$. These and the construction of $\sigma$ imply 
$f(\a_i) - \sigma(\a_i)= d = \sigma(y)-f(y)$. We also observe that $\forall z \in \Delta$, 
$|f(z)-\sigma(z)| \leq \sigma(y)-f(y) = d$, as $f(z)$ is between $\pi$ and $\rho$ whose midplane is $\sigma$.
Now with the above in hand we are able to argue as follows. 
Assume $\mu$ is a best affine approximating plane. Then first $\mu(\a_i) \geq \sigma(\a_i),~~ i \in [1 \dots k+1]$ 
otherwise the maximal distance increases. On the other hand, writing $y=\sum_{i=1}^{k+1}\gamma_i \a_i$ for constants $\gamma_i\geq 0: \Sigma_{i=1}^{k+1}\gamma_i =1$, we have 
$\mu(y)= \sum_{i=1}^{k+1}\gamma_i \mu(\a_i) \leq \sigma (y)$ otherwise the maximal distance increases.
We deduce that $\mu=\sigma$. This concludes the proof of Theorem~\ref{t:simplexD} in case $f$ is convex.

In the case that $f$ is concave, since $\pi = {\rm AS}((a_1, f(\a_1)) \dots (\a_{k+1},f(\a_{k+1})))$ and the graph of $f$ over $\Delta$ lies on or above $\pi$, we can repeat the convex argument replacing $d\text{ by }-d$, $\inf\text{ by }\sup$ and interchanging $\le\text{ and } \ge$. 
\begin{comment}
apply the proof to 
\end{comment} 
Alternatively note $-f$ is convex and one checks that the negative of the solution for $-f$ is just $\frac{\pi+\rho}{2}$. 
\begin{comment}
concave $f$
\end{comment} 
\end{proof}
Note that even if $f$ is not convex (or concave) and not even smooth the method produces a best uniform approximation for the $k+2$ points comprised of the $k+1$ points $(\a_1,f(\a_1))\ldots,(\a_{k+1},f(\a_{k+1}))$ and $(\y,f(\y))$.

\begin{figure}
\centering
\includegraphics[width=0.5\textwidth]{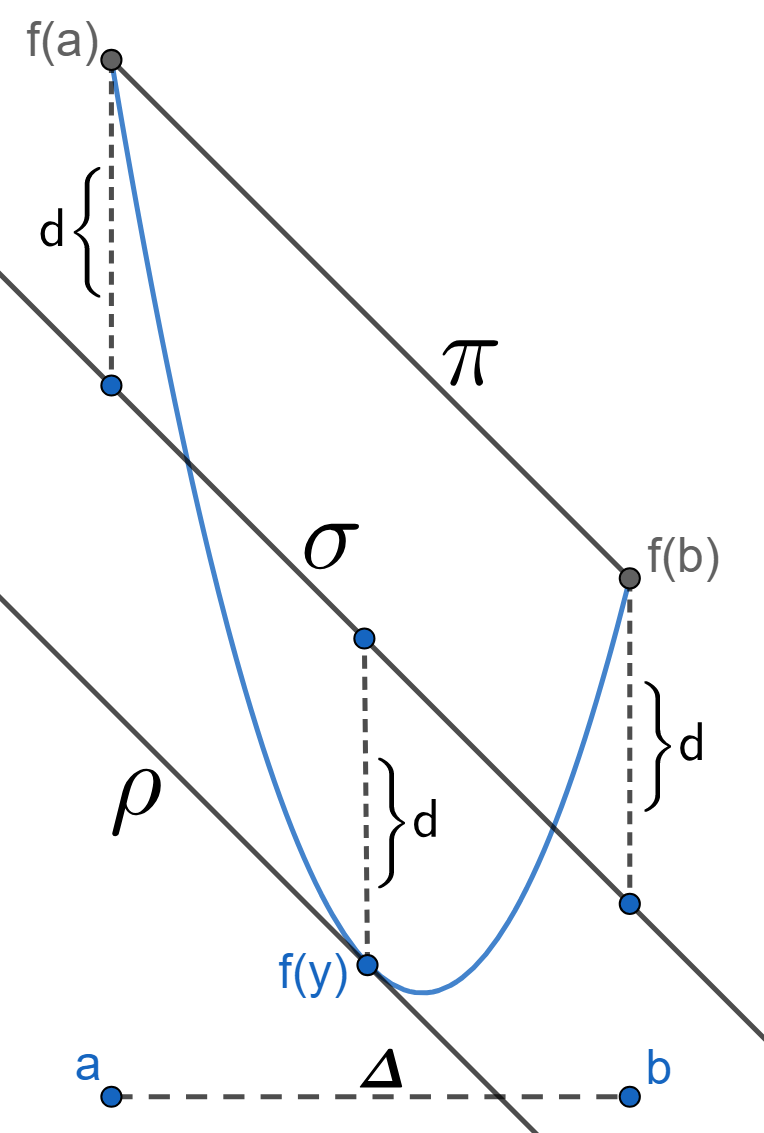}
\caption{Illustration of Theorem~\ref{t:simplexD}. The 
secant $\pi$ is the convex hull of $k+1$ values $(\x_i,f(\x_i))$ and $\rho$ the supporting hyperplane. The convexity of $f$ leads to the fact that the plane, $\sigma$, with same slope halfway between $\pi$ and $\rho$ is the best affine approximation.}
\end{figure}

\section{Theorem~\ref{t:simplexD}: The case $k=1$.}
\setcounter{equation}{0}
\setcounter{theorem}{0}

For $k=1$, our main result Theorem \ref{t:simplexD} provides an extension of the classical Chebyshev equioscillation theorem, see Theorem~\ref{t:equioc1} below, for linear approximants with an explicit unique formula for the best approximant. 
\medskip

We now explain this. 

\subsection{Chebyshev systems and Chebyshev-Markov equioscillation}

We will work with the real interval $[p,q]$ where $p<q$ and the space of real valued continuous functions $f:[p,q]\to \mathbb R$. As per convention, we denote this space by $C([p,q])$.

A set of $l+1,\, l\geq 0$ functions $\left\{u_j(t)\right\}_{j=0}^{l}$, $u_j\in C([p,q])$ is called a Chebyshev (Haar) system on the interval $[p,q]$ if any linear combination
\[
u(t)=\sum_{j=0}^{l}c_ju_j(t),\, t\in [p,q]
\] 
with not all coefficients $c_j$ zero, has at most $l$ distinct zeros in $[p,q]$. 
\footnote{Alfred Haar, 1885-1933, was a Hungarian mathematician. In 1904 he began to study at the University of Göttingen. His doctorate was supervised by David Hilbert.}

The $l+1$ dimensional subspace $U_l\subset C([p,q)])$ spanned by a Chebyshev system 
$\left\{u_j(t)\right\}_{j=0}^l$ on $[p,q]$ 
defined by
\[
U_l:=\left\{u(t):\, u(t)=\sum_{j=0}^{l}c_ju_j(t)\right\}
\]
is called a Chebyshev space on $[p,q]$. 
\medskip

The classical Chebyshev equioscillation theorem, see for example [\cite{Kr}, Chapter 9, Theorem 4.4] and \cite{Sch} is the following:

\begin{theorem}
Let $f\in C([p,q])$ and let $U_l$ be a $l+1$ Chebyschev space on an interval $[p,q]$. Then $\hat{u}\in U_l$ is a best uniform approximant to $f$ on $[p,q]$, that is 
$u$ satisfies
\[
\left\|f-\hat{u}\right\|_{\infty [p,q]}={\inf}_{u\in {U}_l}\left\|f-u\right\|_{\infty [p,q]}
\] 
if and only there exists $l+2$ points $\left\{x_1,..., x_{l+2}\right\}$
with $p\leq x_1<...<x_{l+2}\leq q$ such that 
\beq
f(x_i)-\hat{u}(x_i)=w(-1)^i\|f-\hat{u}\|_{\infty [p,q]},\, w=\pm 1.
\eeq
\label{t:equioc1}
\end{theorem}

Motivated by Theorem~\ref{t:equioc1} % 4.1 
we have:
%\ref{t:equioc1} something odd. 

\subsection{Chebyshev equioscillation}

Let $f\in C([p,q)]$. A Chebyshev polynomial $\hat{h}_l$ (when it exists) of degree at least $l\geq 1$, is the polynomial which best uniformly approximates $f$ on $[p,q]$, i.e, 
\newcommand{\argmin}[1]{\underset{#1}{\operatorname{arg}\,\operatorname{min}}\;} 
2

\[
\hat{h}_l=\argmin{h_l\in {\Pi}_l}\left\|f-h_l\right\|_{\infty [p,q]}
\] 
where ${\Pi}_l$ is the set of polynomials of degree at most $l$. 
\medskip

The following is often called the Chebyshev equioscillation theorem.

\begin{theorem}
Let $f\in C([a,b])$. Then 
$\hat{h_l}$
exists if there exists $l+2$ points $\left\{x_1,..., x_{l+2}\right\}$
with $p\leq x_1<...<x_{l+2}\leq q$ such that 
\beq
f(x_i)-\hat{h_l}(x_i)=w(-1)^i\|f-\hat{h_l}\|_{\infty [p,q]},\, w=\pm 1.
\eeq
\label{t:equioc}
\end{theorem}

\subsection{The case $l=1$ of Theorem~\ref{t:equioc}} %4.2}

We now show how Theorem~\ref{t:equioc} %4.2 
with $l=1$ corresponds to Theorem 2.1 for $k=1$.
\medskip

We may assume that $p<q$ and for the moment we do not assume anything about $f:[p,q]\to \mathbb R$.
Define now a linear function $L:[p,q]\to \mathbb R$ from $f$ in terms of parameters $d$ and $y$ to be determined later as follows: 
\beq
L(x)=f(y)+d+m(x-y),\, x\in [p,q].
\label{e:L}
\eeq
Now since $L(y)-f(y)=d$, if we assume $L(p)-f(p) = -d = L(q)-f(q)$ then
\[
d=f(p)-L(p)=f(p)-f(y)-d-m(p-y).
\]
Thus 
\[
d=\frac{f(p)-f(y)}{2}-\frac{m}{2}(p-y).
\]
Also, 
\[
d=f(q)-L(q)=f(q)-f(y)-d-m(q-y).
\]
So: 
\[
\frac{f(p)-f(y)}{2}-\frac{m}{2}(p-y)=\frac{f(q)-f(y)}{2}-\frac{m}{2}(q-y).
\]
Thus $m=\frac{f(q)-f(p)}{p-q}$ in the case when $f$ or $-f$ is a convex and differentiable function, by the mean value theorem.
\medskip

It is clear that we have proved Theorem~\ref{t:simplexD} once we are able to choose $r$ to maximize $d$ if this is possible for the given $f$. In the case when $f$ is convex, we see that $m=f'(r)$.
It is clear that the argument works when $f$ is concave, which implies $-f$ is convex, 
\medskip

In the case of $k=1$, we could write
\[
\frac{ax+b}{cx+d}==\frac{\left(\left(\frac{a}{c}\right)\left(\x+\frac{d}{c}-\frac{d}{c}\right)+\frac{b}{c}\right)}{\left(\x+\frac{d}c\right)}
\]
\[
=\frac{a}{c}+\frac{\left(\frac{b}{c}-\left(\frac{a}{c}\right)\left(\frac{d}{c}\right)\right)}{\left(\x+\frac{d}{c}\right)}
\]
which is either an upward or downward hyperbola for which the secant between $[p,q]$ lies above (or below) the curve. Then following the "mean value " argument, leads to $r$ with derivative =slope and via our visualization leads to a line with the same slope but halfway between the secant and tangent to $r$.

\section{Connections of Theorem~\ref{t:simplexD} to graphics}
\setcounter{equation}{0}

One consequence of Theorem~\ref{t:simplexD} gives an interesting connection to graphics. We provide our ideas below.
\medskip

We are given a flat object $\mathcal O$ 
and want to render it from a given perspective, camera setup.
The resulting image is 
a projective (AKA perspective or homography) transformation of $\mathcal O$, $P(\mathcal O)$.
Thus, each pixel (color) in $\mathcal O$, (an ordered pair $(x_1,y_1)$ in $\mathbb R^2$) is transformed to a new pixel location (point) 2 via the action of $P$
\begin{equation}
P:\left[    
\begin{matrix} 
x_1
\\y_1 
\end{matrix}\right] \mapsto
\frac{1}{d x_1+e y_1 +j}
\begin{bmatrix}
{a_1 x_1+b_1 y_1 +c_1}\\
{a_2 x_1+b_2 y_1 +c_2}
\end{bmatrix}
\label{e:projtran}.
\end{equation}

Practically, to render this object which can have 10's of millions of pixels it is useful to have a good fast approximation of the transformation not entailing division which can cause numerical instabilities. 
Thus, we seek affine approximants to $P$. One known method is to simply take the affine approximant to be the linear terms of the Taylor expansion around one point, the tangent approximation.

We aim to provide a better 2d-affine approximant to $P$ than the tangent approximation. From Equation~\ref{e:projtran} the components of $P(\x)$ are given by
\begin{equation} 
f_i(\x) := \frac{a_i x + b_i y + c_i}{d x + e y + j}, \,\,\,\, i = 1, 2 .
\label{e:Pcomp}
\end{equation}

From our main Theorem~\ref{t:simplexD} we know that if we can find a triangle, $\Delta$, containing (a large part of) $\mathcal O$ such that each $f_i(\x)$ 
is convex or concave then there exist  Min-Max affine approximations ${{\bm \alpha}_i}^T \x + \beta_i$ to each $f_i$. Then forming the affine transformation $A$ of 2-space given by 

\begin{equation}
A(\x):= 
\begin{bmatrix} 
{{\bm \alpha}_1}^T \x + \beta_1 \\ 
{{\bm \alpha}_2}^T \x + \beta_2 
\end{bmatrix},
\label{e:twodaffine}
\end{equation}
provides, component wise, the best uniform affine approximants on $\Delta$.

Since a differentiable function is convex or concave on a domain exactly when it{'}s Hessian is respectively positive or negative semidefinite the following is a key tool for the application of our main result.

\begin{comment}
The key to characterizing the regions in $\RR^2$ over which $P$ is convex or concave is the following fact:
\end{comment}
\begin{prp}
\label{p:twodposdef}
A symmetric 
$2 x 2$ matrix 
$D = 
\begin{bmatrix} 
p_{11} & p_{12} \\ 
p_{12} & p_{22} 
\end{bmatrix}$
is \textnormal{ positive }\textnormal{or negative}   semidefinite  according as 
$signum(p_{11}) = {\pm} \textnormal{signum }(p_{11} p_{22} - p_{12}^2).$
\end{prp}

The Hessian of 
\begin{equation}
   \frac{\alpha X+\beta Y+\gamma}{X+\delta} = \alpha + \frac{\beta Y - \alpha \delta + \gamma }{X+\delta}
\label{frac}
\end{equation}
is
$$
H=\frac{1}{(X+\delta)^3}
\begin{bmatrix} 
2 ( \beta Y-\alpha \delta  + \gamma) & -\beta(X + \delta)\\
-\beta(X + \delta) & 0
\end{bmatrix} 
$$
with determinant 
$\frac{-\beta^2}{(X + \delta)^4}$.
The sign of 
$H[1,1]=2 ( \beta Y -\alpha \delta + \gamma)$  depends on $Y \begin{matrix} < \\ > \end{matrix} \frac{\alpha \delta-\gamma}{\beta}$.

When we combine Proposition~\ref{p:twodposdef} with the 
Hessian of Eq.~\ref{frac} we see that $\frac{\alpha X+\beta Y+\gamma}{X+\delta}$ in each of the connected components of $\RR^2\setminus \{ \{X=-\delta\} \cup \{Y=\frac{\alpha \delta-\gamma}{\beta}\}\}$  is either  convex or  concave.
With the 2 coordinate functions of the projective transformation ($X$ and $Y$) we end with 6 regions in 
$\RR^2$ where Theorem~\ref{t:simplexD} holds, Figure~\ref{fig:M1}.

\begin{figure}
\centering
\scalebox{.7}{
\begin{tikzpicture}
 \draw [very thick,color = blue, dashed] (-1.5,-0.5) -- (9,-0.5);
     \draw[fill=red!10] (0,4) -- (0.4,0) -- (2,4.3) -- cycle;
     \draw[fill=yellow!10] (1.25,-3.5) -- (1.75,-1.4) -- (2.5,-3.7) -- cycle;
    \draw[fill=blue!10] (5,0) -- (8,1) -- (7.7,3) -- cycle;
\draw [very thick](3,-4) -- (3,5);
\draw [very thick, color = red, dashed](-1.5,3.35) -- (9,3.35);
       \end{tikzpicture}
}
\caption{The black line is where the projective transformation is infinite, $X= -\delta$. The blue line is where the 
$\frac{\partial{}}{\partial{x^2}}$ term of the  Hessian of the $X$ transformation is $0$ and the red line is where the $\frac{\partial{}}{\partial{x^2}}$ term of the  Hessian
of the $Y$ transformation is $0$.
The blue and yellow simplices define domains of transformations that are convex (or concave) while the red simplex is not.} 
\label{fig:M1}
\end{figure}
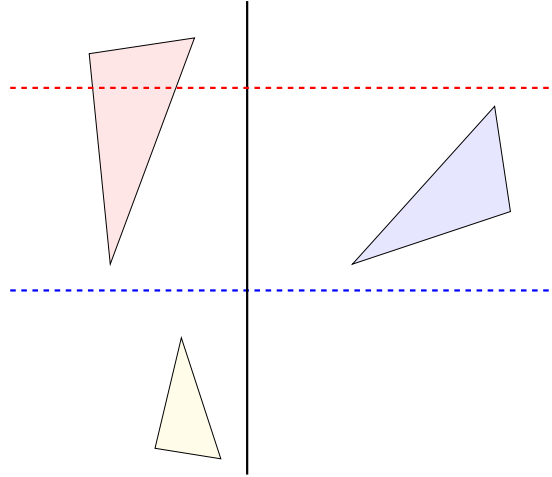

In order to transform a general projective transformation as in
Eq. ~\ref{e:twodaffine}
to the form in Eq.~\ref{frac}, with a simple denominator,
we rotate the axes so that $X$ is in the direction of 
$\begin{bmatrix}d\\
e
\end{bmatrix}
$ and normalize so that the coefficient of $X$ in the denominator is 1.
\medskip

In Figure~\ref{teddy} we show a few examples of the various affine approximations $A$ of projective transformations $P$ using this idea.
\medskip
 
The first column is the original image $\mathcal O$,
the second column is the image transformed by a projective transformation $P$, modelling a new viewpoint, the third column is the transformation based on the affine approximation $A$  of $P$  using Theorem~\ref{t:simplexD} and a user defined triangle while the fourth column is computed using a Taylor series approximation of $P$ around the image's center. The images in column 3 
should visually look closer to those in column 2 than those in column 4, the Taylor approximation.

\begin{figure}
\centering

\includegraphics[width=0.24\textwidth]{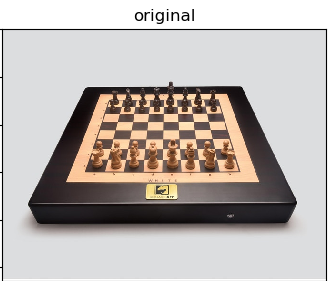}
\includegraphics[width=0.24\textwidth]{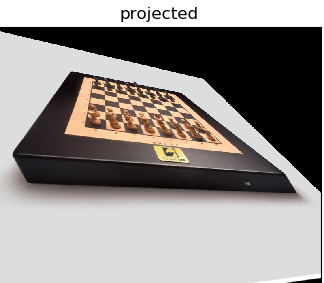}
\includegraphics[width=0.24\textwidth]{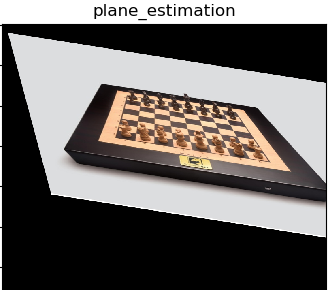}
\includegraphics[width=0.24\textwidth]{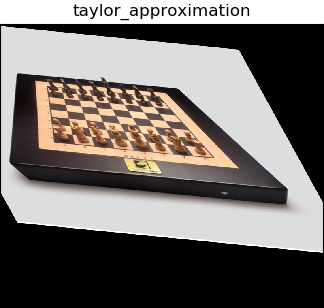}\\
\includegraphics[width=0.24\textwidth]{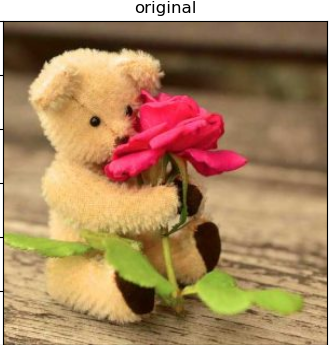}
\includegraphics[width=0.24\textwidth]{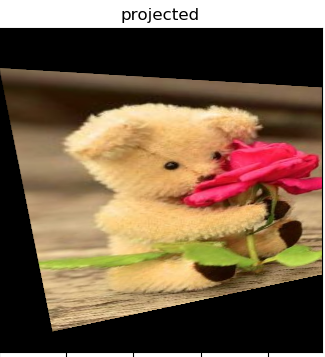}
\includegraphics[width=0.24\textwidth]{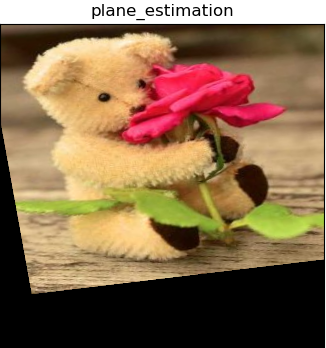}
\includegraphics[width=0.24\textwidth]{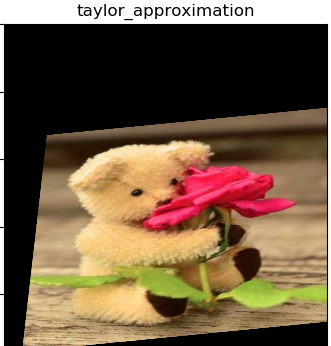}\\
\includegraphics[width=0.24\textwidth]{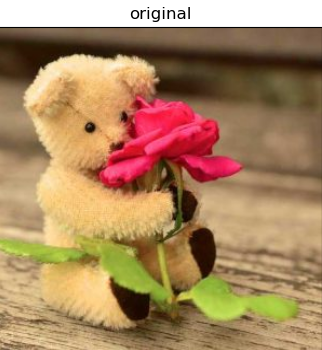}
\includegraphics[width=0.24\textwidth]{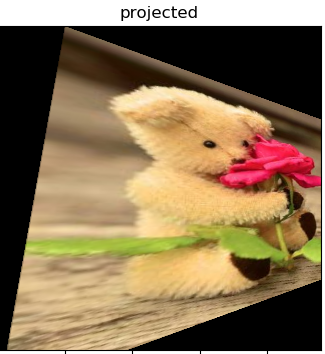}
\includegraphics[width=0.24\textwidth]{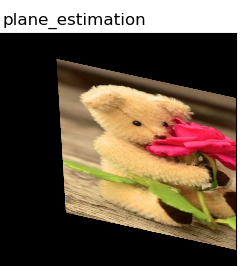}
\includegraphics[width=0.24\textwidth]{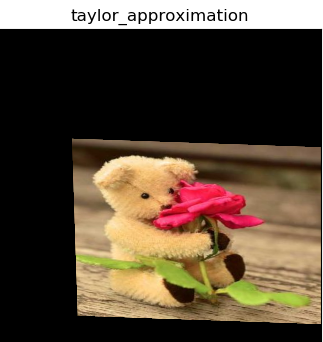}\\
\includegraphics[width=0.24\textwidth]{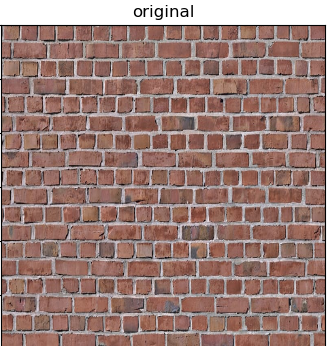}
\includegraphics[width=0.24\textwidth]{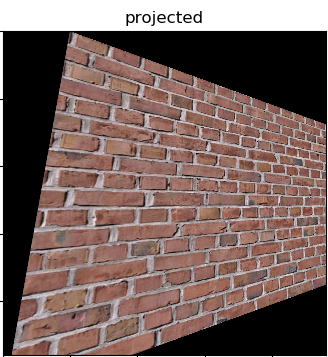}
\includegraphics[width=0.24\textwidth]{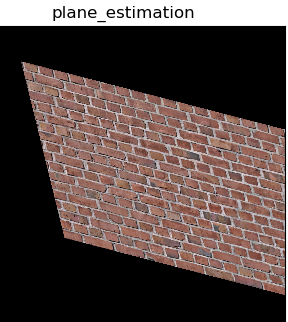}
\includegraphics[width=0.24\textwidth]{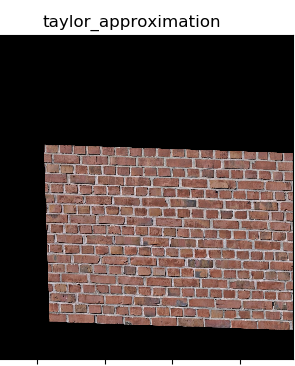}
\caption{Original image / projective transformation of the image/ best uniform affine approximation/ Taylor expansion.}
\label{teddy}
\end{figure}

\pagebreak

\end{document}